\documentclass[11pt,english]{article}

\usepackage[T1]{fontenc}
\usepackage[utf8]{inputenc}
\usepackage[margin=1in]{geometry}
\usepackage{babel}
\usepackage{amsmath,amssymb,amsfonts,amsthm,mathrsfs,mathtools}
\usepackage[dvipsnames]{xcolor}
\usepackage{enumitem}
\usepackage{microtype}
\usepackage{hyperref}
\usepackage[nameinlink,noabbrev]{cleveref}
\usepackage{currfile}

\hypersetup{
  colorlinks=true,
  linkcolor=MidnightBlue,
  citecolor=ForestGreen,
  urlcolor=MidnightBlue,
  pdftitle={The Minimal Abstract Robust Subdifferential},
  pdfauthor={M.D. Voisei}
}

\newcommand{\Rbar}{\overline{\mathbb{R}}}
\newcommand{\Graph}{\operatorname{Graph}}
\newcommand{\cl}{\operatorname{cl}}
\newcommand{\co}{\operatorname{co}}
\newcommand{\fen}{\partial_{\!F}}
\newcommand{\efen}{\partial_{\!F}^{\varepsilon}}
\newcommand{\fsub}{\widehat{\partial}}
\newcommand{\lsub}{\partial_{\!L}}
\newcommand{\csub}{\partial_{\!C}}
\newcommand{\gsub}{\partial_{\!G}}
\newcommand{\ind}{\iota}
\newcommand{\wstar}{w^{*}}
\DeclareMathOperator{\dist}{dist}

\DeclareMathOperator*{\Limsup}{Limsup}

\theoremstyle{plain}
\newtheorem{theorem}{Theorem}
\newtheorem{proposition}[theorem]{Proposition}
\newtheorem{lemma}[theorem]{Lemma}
\newtheorem{corollary}[theorem]{Corollary}

\theoremstyle{definition}
\newtheorem{definition}[theorem]{Definition}
\newtheorem{example}[theorem]{Example}

\theoremstyle{remark}
\newtheorem{remark}[theorem]{Remark}

\Crefname{theorem}{theorem}{theorems}
\Crefname{theorem}{Theorem}{Theorems}
\Crefname{proposition}{proposition}{propositions}
\Crefname{proposition}{Proposition}{Propositions}
\Crefname{lemma}{lemma}{lemmas}
\Crefname{lemma}{Lemma}{Lemmas}
\Crefname{corollary}{corollary}{corollaries}
\Crefname{corollary}{Corollary}{Corollaries}
\Crefname{definition}{definition}{definitions}
\Crefname{definition}{Definition}{Definitions}
\Crefname{example}{example}{examples}
\Crefname{example}{Example}{Examples}
\Crefname{remark}{remark}{remarks}
\Crefname{remark}{Remark}{Remarks}

\begin{document}

\title{The Minimal Robust Core of Abstract Subdifferentials}
\author{M.D. Voisei}
\date{}
\maketitle

\begin{abstract}
This paper introduces the metric-dependent core subdifferential, a local
robust affine-support construction for extended-real functions on metrizable
topological vector spaces.  A core subgradient is a continuous linear slope for
which the affine lower support holds on metric balls with an error negligible
relative to the supporting radius, along arbitrarily fine admissible scales.
The main minimality results show that these slopes are unavoidable: on complete
metrizable spaces they belong to the graph closure of any subdifferential
satisfying a local-minimum principle together with a mild stability condition
under metric-distance perturbations, and on Banach spaces they belong to the
nearby graph closure of any abstract subdifferential satisfying the usual fuzzy
minimum principle.  For the norm metric, the construction contains the
Fr\'echet subdifferential, coincides with the Fenchel subdifferential on convex
functions, is contained in the limiting subdifferential whenever the relevant
Fr\'echet fuzzy calculus is available, and is contained in the
Clarke--Rockafellar subdifferential for lower semicontinuous functions on
Banach spaces.  The paper also records metric-dependence phenomena, strict
comparison examples, constrained optimality and variational-inequality
conditions, a scale-slope/error-bound characterization, and the relation with
Goldstein-type stationarity used in finite-dimensional nonsmooth optimization.
\end{abstract}



\section{Introduction}

Generalized differentiation provides several non-equivalent ways to encode
first-order information for nonsmooth and extended-real-valued functions.
Convex analysis uses the Fenchel subdifferential; variational analysis uses
Fr\'echet, limiting, and Clarke--Rockafellar constructions; and algorithmic
nonsmooth optimization often works with relaxed neighborhood objects such as
Goldstein subdifferentials.  These constructions differ in their support
properties, calculus rules, closure behavior, and stationarity interpretation;
see, for example,
\cite{MR709590,RockafellarWets1998,Mordukhovich2006,Ioffe2017}.

The aim of this paper is to isolate a small local-support object that is forced
by weak and standard subdifferential principles.  Given a compatible metric
\(d\), a continuous linear functional \(x^{*}\) is declared to be a core
subgradient of \(f\) at \(x\) when the affine support
\[
 f(x)+\langle u-x,x^{*}\rangle\le f(u)+\varepsilon
\]
holds on a \(d\)-ball around \(x\) whose radius is large compared with
\(\varepsilon\), along errors \(\varepsilon\downarrow0\).  Thus the core
condition is stronger than merely having approximate support at one point, but
weaker than exact local support on a fixed neighborhood.  It records precisely
those slopes for which the affine perturbation \(f-x^{*}\) has sublinear best
local descent along arbitrarily small scales.

The construction is intentionally metric-dependent.  Compatible metrics that
are locally bi-Lipschitz equivalent give the same core, but arbitrary compatible
metrics need not do so.  In normed spaces all normed-space comparisons in this
paper use the norm metric; equivalent norms therefore give the same object.  In
the general metrizable-vector-space results, no translation invariance of the
chosen metric is assumed unless explicitly stated.

The main minimality principle is proved in two forms.  First,
\Cref{thm:min-sum} treats complete metrizable topological vector spaces and
subdifferentials satisfying a Fermat rule plus stability under perturbations by
continuous linear functions and the metric distance.  Second,
\Cref{thm:min-abstract} gives the Banach-space version for abstract
subdifferentials satisfying the usual fuzzy nearby-point minimum principle.  In
both cases, every core subgradient belongs to the corresponding graph outer
limit of the chosen subdifferential.  Nets are used in the definitions because
the relevant graph closures need not be first countable; in first-countable
situations they can be replaced by sequences.

The normed-space comparisons then locate the core among standard generalized
gradients.  The Fr\'echet subdifferential is always contained in the core.  When
the Fr\'echet subdifferential satisfies the fuzzy minimum principle, for example
for proper lower semicontinuous functions on Asplund Banach spaces, the core is
contained in the limiting subdifferential.  For convex functions the core
coincides with the Fenchel subdifferential.  For lower semicontinuous functions
on Banach spaces it is contained in the Clarke--Rockafellar subdifferential.
Examples show that the inclusions can be strict: the core excludes slopes that
come only from convexification, neighborhood relaxation, or sharp local maxima.

The rest of the paper develops consequences of these comparisons.  The
\(\varepsilon\)-minimality-radius and scale-slope characterizations give an
error-bound interpretation of the definition.  The core normal cone gives direct
necessary conditions for constrained minimization and local variational
inequalities.  Finally, in finite dimensions the Clarke comparison implies that
core stationarity is stronger than Goldstein stationarity at every fixed radius,
while simple cusp examples show that the converse fails.

\section{The core subdifferential and abstract minimality}

Let $(X,s)$ be a real metrizable topological vector space with topological
dual space $X^{*}$ and duality product
$\langle x,x^{*}\rangle:=x^{*}(x)$, $x\in X$, $x^{*}\in X^{*}$.  Let
$d$ be a metric on $X$ compatible with $s$, meaning that $d$ generates the
linear topology $s$ of $X$.  Denote the closed $d$-ball by
\[
\overline{B}_d(x;r):=\{y\in X\mid d(x,y)\le r\},\qquad x\in X,\ r\ge0.
\]
For a multifunction $T:A\rightrightarrows B$, write
$\Graph(T):=\{(a,b)\in A\times B\mid b\in T(a)\}$. 

\begin{definition}[core subdifferential]
For $f:X\to\Rbar$ and $x_{0}\in X$, the \emph{$d$-core subdifferential}
$\eth_d f(x_{0})\subset X^{*}$ is defined as follows: $x_{0}^{*}\in\eth_d f(x_{0})$
if $f(x_{0})\in\mathbb{R}$ and there exists
$(\delta_{\varepsilon})_{\varepsilon>0}\subset(0,+\infty)$ such that
\[
\limsup_{\varepsilon\downarrow0}\delta_{\varepsilon}/\varepsilon=+
\infty
\]
and, for every $\varepsilon>0$ and every $x\in \overline{B}_d(x_{0};\delta_{\varepsilon})$,
\begin{equation}\label{core}
f(x_{0})+\langle x-x_{0},x_{0}^{*}\rangle\le f(x)+\varepsilon.
\end{equation}
If $f(x_{0})\notin\mathbb{R}$, set $\eth_d f(x_{0})=\emptyset$.
\end{definition}

The definition would be unchanged if open balls were used instead of closed balls, since one may always slightly shrink the radius. We use closed balls because this form is more convenient for the arguments that follow, especially for the application of Ekeland's variational principle. 

\begin{remark}[dependence on the metric]\label{rem:metric-dependence}
The metric $d$ is part of the definition. 

The construction is unchanged if \(d\) is replaced by a metric \(\rho\)
which is locally bi-Lipschitz equivalent to \(d\) near the diagonal.
More precisely, assume that for every \(x_{0}\in X\) there exist a
neighborhood \(U\) of \(x_{0}\) and constants \(a,b>0\) such that
\[
a\,d(x,y)\leq \rho(x,y)\leq b\,d(x,y),
\qquad x,y\in U .
\]
Then, for all sufficiently small \(r>0\), the corresponding balls satisfy
\[
\overline{B}_d(x_{0};r/b)\subset \overline{B}_{\rho}(x_{0};r)
\subset \overline{B}_d(x_{0};r/a).
\]
Hence the defining neighborhoods in the definition of
\(\eth_{d}f(x_{0})\) and \(\eth_{\rho}f(x_{0})\) differ only by fixed
multiplicative changes of radius. Such changes do not affect the condition
\(
\limsup_{\varepsilon\downarrow 0}
\delta_{\varepsilon}/{\varepsilon}=+\infty
\). 
Consequently, 
\[
\eth_{d}f(x_{0})=\eth_{\rho}f(x_{0}).
\]

The $d$-core subdifferential is not invariant under arbitrary changes of topologically compatible metric. For example, let $X=\mathbb R$, $f\equiv 0$, and $x_0=0$. For the usual metric $d_E(x,y)=|x-y|$, the functional $x^*=1$ does not belong to $\eth_{d_E}f(0)$. 

Indeed, the defining inequality would require, for each sufficiently small $\varepsilon>0$, the existence of a radius $\delta_\varepsilon>0$ such that
\(x\leq \varepsilon\) for all $x\in \overline{B}_{d_0}(0;\delta_\varepsilon)=[-\delta_\varepsilon,\delta_\varepsilon]$. 

Taking \(x=\delta_\varepsilon\) gives \(\delta_\varepsilon\leq \varepsilon\). Hence \(\delta_\varepsilon/{\varepsilon}\leq 1\), which is incompatible with the core condition

\[
\limsup_{\varepsilon\downarrow 0}\delta_\varepsilon/{\varepsilon}=+\infty.
\]

Thus \(1\notin \eth_{d_E}f(0)\).

On the other hand, consider the compatible metric
\(d(x,y)=|x-y|^{1/2}\). 
Then 
\(
\overline{B}_d(0;\delta)=[-\delta^2,\delta^2]\). 
Choosing
\(
\delta_\varepsilon=\sqrt{\varepsilon}
\)
gives
\(
\overline{B}_d(0;\delta_\varepsilon)=[-\varepsilon,\varepsilon].
\)
Thus, for every \(x\in \overline{B}_d(0;\delta_\varepsilon)\), we have
\(
x\leq \varepsilon.
\)
Moreover, \(\delta_\varepsilon/{\varepsilon}=\varepsilon^{-1/2}\to+\infty\) as $\varepsilon\downarrow0$. 
Therefore $ 1\in \eth_d f(0)$. 

This shows that the construction depends on the metric and is not invariant under arbitrary compatible changes of metric.

When \(X\) is a normed space, we write \(\eth_{\|\cdot\|}\) simply as \(\eth\). This causes no ambiguity for equivalent norms: if \(\|\cdot\|_1\) and \(\|\cdot\|_2\) are equivalent norms, then their induced metrics are globally bi-Lipschitz equivalent, and therefore they generate the same core subdifferential. Consequently, all statements below, e.g. comparing \(\eth\) with the Fenchel subdifferential, are made with respect to the norm metric.
\end{remark}

By a \emph{subdifferential} we mean any map
$\partial:\Rbar^{X}\times X\rightrightarrows X^{*}$.  We write
$x_{i}\to^{s}x$ for convergence in $(X,s)$.

If $\tau$ is a topology on $X^{*}$, define
\begin{equation}\label{eq:Limsup}
x^{*}\in s\times\tau-\Limsup_{y\to x}\partial f(y)
\end{equation}
when there is a net $((y_{i},y_{i}^{*}))_{i}\subset\Graph(\partial f)$ such
that $y_{i}\to^{s}x$, $y_{i}^{*}\to^{\tau}x^{*}$, and
$f(y_{i})\to f(x)$.  If the product topology is first countable on the relevant
sets, this net can be replaced by a sequence. 
In particular, $\partial f(x)\subset s\times\tau-\Limsup_{y\to x}\partial f(y)$, $x\in X$. 

We identify a functional $x^{*}\in X^{*}$ with the continuous linear
function $u\mapsto\langle u,x^{*}\rangle$ when writing sums such as
$f+x^{*}$.

Let $\mu^{*}$ be a linear topology on $X^{*}$. For $x\in X$, define
$S_{\mu^{*}}(\partial,x)$ as the class of all functions $f:X\to\overline{\mathbb R}$
for which there exist a $\mu^{*}$-bounded set $M_{f,x}\subset X^{*}$ and a
neighborhood $V_{f,x}$ of $x$ such that $f$ is lower semicontinuous (l.s.c.) on
$V_{f,x}$ and, for all $y\in V_{f,x}$, $\lambda\ge0$, and $x^{*}\in X^{*}$,
\begin{equation}\label{eq:S-gen}
	\partial\bigl(f+x^{*}+\lambda d(\cdot,y)\bigr)(y)
	\subset \partial f(y)+x^{*}+\lambda M_{f,x}.
\end{equation}

\begin{theorem}\label{thm:min-sum}
Let $(X,s)$ be a metrizable topological vector space and let $d$ be a
compatible metric on $X$ such that $(X,d)$ is complete.  Let $\mu^{*}$ be a
linear topology on $X^{*}$ and let
$\partial:\Rbar^{X}\times X\rightrightarrows X^{*}$ be a subdifferential with
property
\begin{description}[leftmargin=2.8em,style=nextline]
\item[{\rm[M]}] If $h:X\to\Rbar$, $x_{0}$ is a local minimum point of $h$,
and $h(x_{0})\in\mathbb{R}$, then $0\in\partial h(x_{0})$.
\end{description}
Then, for every $x\in X$ and every $f\in S_{\mu^{*}}(\partial,x)$,
\[
\eth_d f(x)\subset s\times\mu^{*}-\Limsup_{y\to x}\partial f(y).
\]
\end{theorem}

\begin{proof}
Let $x_{0}\in X$, $f\in S_{\mu^{*}}(\partial,x_{0})$, and
$x_{0}^{*}\in\eth_d f(x_{0})$.  Choose sequences
$\varepsilon_{n}>0$, $\delta_{n}>0$ such that
\[
\varepsilon_{n}\to0,\qquad \varepsilon_{n}/{\delta_{n}}\to0,
\]
and
\begin{equation}\label{eq:L}
f(x_{0})+\langle x-x_{0},x_{0}^{*}\rangle\le f(x)+\varepsilon_{n}
\end{equation}
whenever $x\in \overline{B}_d(x_{0};\delta_{n})$.

If $\inf_{n}\delta_{n}>0$, then $x_{0}$ is a local minimum point of
$F:=f-x_{0}^{*}$: the inequalities above hold on a fixed ball after passing to
all sufficiently large $n$, and then $\varepsilon_n\downarrow0$ along a
subsequence.  By [M], $0\in\partial(f-x_{0}^{*})(x_{0})$.  Applying
\eqref{eq:S-gen} with $\lambda=0$ and perturbation $-x_{0}^{*}$ gives
$0\in\partial f(x_{0})-x_{0}^{*}$, hence
$x_{0}^{*}\in\partial f(x_{0})\subset s\times\mu^*-\Limsup_{y\to x_0}\partial f(y)$.

Assume now, after passing to a subsequence, that $\delta_{n}\to0$.  Choose
$\delta>0$ so that $\overline{B}_d(x_{0};\delta)$ is contained in the neighborhood
where $f$ is l.s.c. and \eqref{eq:S-gen} holds, and assume
$\delta_{n}\le\delta$.  Set $F:=f-x_{0}^{*}$.  From \eqref{eq:L},
\[
-\infty<F(x_{0})\le \inf_{\overline{B}_d(x_{0};\delta_{n})}F+\varepsilon_{n}.
\]
Ekeland's variational principle \cite{Ekeland1974}, applied to the complete
metric space $\overline{B}_d(x_{0};\delta_{n})$, yields
$y_{n}\in \overline{B}_d(x_{0};\delta_{n})$ such that
$F(y_{n})\le F(x_{0})$, $d(y_{n},x_{0})\le\delta_{n}/2$, and $y_n$ minimizes
$F+(2\varepsilon_n/\delta_n)d(\cdot,y_n)$ locally in $X$.  The last assertion
follows because the ball $B_d(y_n;\delta_n/2)$ is contained in
$\overline B_d(x_0;\delta_n)$.

Property [M] and the defining inclusion \eqref{eq:S-gen} applied with
$x^*=-x_0^*$, $\lambda=2\varepsilon_{n}/\delta_{n}$, and $y=y_n$ give
\[
0\in\partial\biggl(f-x_{0}^{*}+(2\varepsilon_{n}/\delta_{n})d(\cdot,y_{n})\biggr)(y_{n})
\subset
\partial f(y_{n})-x_{0}^{*}+(2\varepsilon_{n}/\delta_{n})M,
\]
for some $\mu^*$-bounded $M=M_{f,x_0}$.

Thus there are $y_{n}^{*}\in\partial f(y_{n})$ and $m_{n}^{*}\in M$ such that
$y_{n}^{*}=x_{0}^{*}-(2\varepsilon_{n}/\delta_{n})m_{n}^{*}$.  Since $M$ is
$\mu^{*}$-bounded and $2\varepsilon_{n}/\delta_{n}\to0$, we have
$y_{n}^{*}\to^{\mu^{*}}x_{0}^{*}$.  Also $y_{n}\to^{s}x_{0}$.  Finally, $F(y_{n})\le F(x_{0})$ gives
\[
\limsup_n f(y_n)\le f(x_0),
\]
because $\langle y_n-x_0,x_0^*\rangle\to0$, while the lower semicontinuity of
$f$ at $x_0$ gives $\liminf_n f(y_n)\ge f(x_0)$.  Hence
$f(y_{n})\to f(x_{0})$, and therefore
$x_{0}^{*}\in s\times\mu^{*}-\Limsup_{y\to x_{0}}\partial f(y)$.
\end{proof}

\begin{remark}\label{rem:sum-rule-context}
A common way to verify the defining inclusion \eqref{eq:S-gen} is to prove a sum rule of the form
\begin{equation}\label{eq:S}
\partial\bigl(f+x^{*}+\lambda d(\cdot,y)\bigr)(y)
\subset\partial f(y)+x^{*}+\lambda\partial d(\cdot,y)(y),
\end{equation}
where $\partial d(\cdot,y)(y)$ is $\mu^{*}$-bounded.  For example, in a normed space with
$d(x,y)=\|x-y\|$, the set $\partial d(\cdot,y)(y)$ is contained in the dual unit ball 
$B_{X^{*}}$ for subdifferentials that agree with, or are dominated by, the Fenchel subdifferential on continuous convex functions. 
\end{remark}

Let $(X,\|\cdot\|)$ be a Banach space and let
$\emptyset\ne\mathcal{F}\subset\Rbar^{X}$.  Following the terminology of
Thibault \cite{MR1357833} and Z\u{a}linescu \cite{MR1921556}, a multifunction
$\partial:\Rbar^{X}\times X\rightrightarrows X^{*}$ is an \emph{abstract
subdifferential} (or quasi presubdifferential) on $\mathcal{F}$ if it satisfies
\begin{description}[leftmargin=2.8em,style=nextline]
\item[{\rm[P]}] If $f\in\mathcal{F}$, $g:X\to\mathbb{R}$ is convex and
continuous, $x\in X$, $f(x)$ is finite, $f$ is l.s.c. near $x$, and $x$ is a
local minimum point of $f+g$, then
\begin{equation}\label{eq:abs-sub}
0\in s\times w^*-\Limsup_{y\to x}\partial f(y)+\fen g(x),
\end{equation}
where $\fen$ denotes the Fenchel subdifferential, s denotes the strong topology of $X$, and $w^*$ stands for the weak-star topology of $X^{*}$.
\end{description}

In a normed space we write $\overline{B}(x_{0};r):=\{x\in X\mid\|x-x_{0}\|\le r\}$, 
use "$\to$" for norm convergence in $X$ or $X^*$, and denote by $s^*$ the norm topology of $X^*$.

\begin{theorem}\label{thm:min-abstract}
Let $(X,\|\cdot\|)$ be a Banach space and let
$\partial:\Rbar^{X}\times X\rightrightarrows X^{*}$ be an abstract
subdifferential on $\mathcal{F}$.  Then, for every $x\in X$ and every
$f\in\mathcal{F}$ that is l.s.c. near $x$,
\begin{equation}\label{eq:mars2}
\eth f(x)\subset s\times w^*-\Limsup_{y\to x}\partial f(y).
\end{equation}
\end{theorem}

\begin{proof}
Let $x_{0}\in X$, $f\in\mathcal{F}$ be l.s.c. near $x_0$, and $x_{0}^{*}\in\eth f(x_{0})$.  Choose
$\varepsilon_{n}>0$ and $\delta_{n}>0$ such that
$\varepsilon_{n}\to0$, $\varepsilon_{n}/\delta_{n}\to0$, and, for all $x\in \overline{B}(x_{0};\delta_{n})$, 
\begin{equation}\label{eq:L2}
f(x_{0})+\langle x-x_{0},x_{0}^{*}\rangle\le f(x)+\varepsilon_{n}.
\end{equation}

If $\inf_{n}\delta_{n}>0$, then $x_{0}$ is a local minimum point of
$f-x_{0}^{*}$.  Applying \eqref{eq:abs-sub} to the convex continuous function
$g=-x_{0}^{*}$, for which $\fen g(x_0)=\{-x_0^*\}$, gives
$x_{0}^{*}\in s\times w^*-\Limsup_{y\to x_{0}}\partial f(y)$.

Otherwise, pass to a subsequence such that $\delta_{n}\to0$ and such that
$\overline{B}(x_0;\delta_{n})$ is contained in a neighborhood of $x_0$ where $f$ is l.s.c.  Set $F=f-x_{0}^{*}$.  
From the defining inequality,
\[
-\infty<F(x_{0})\le \inf_{\overline{B}(x_{0};\delta_{n})}F+\varepsilon_{n}.
\]

Ekeland's variational principle \cite{Ekeland1974}, applied to the l.s.c. function $F$ on the complete metric space $\overline{B}(x_{0};\delta_{n})$, gives
$y_{n}\in \overline{B}(x_{0};\delta_{n})$ such that $F(y_{n})\le F(x_{0})$,
$\|y_{n}-x_{0}\|\le\delta_{n}/2$, and $y_{n}$ is a local minimum point of
\[
f-x_{0}^{*}+(2\varepsilon_{n}/\delta_{n})\|\cdot-y_{n}\|.
\]
Since $-x_{0}^{*}+(2\varepsilon_{n}/\delta_{n})\|\cdot-y_{n}\|$ is convex and continuous, property [P] gives  
\[
0\in s\times w^*-\Limsup_{z\to y_{n}}\partial f(z)-x_{0}^{*}
+(2\varepsilon_{n}/\delta_{n})\overline{B}_{X^{*}}.
\]
Here $\overline{B}_{X^{*}}=\fen \|\cdot-y_{n}\|(y_n)$ is the closed unit ball
of $X^*$.

Thus there is
$y_{n}^{*}\in s\times w^*-\Limsup_{z\to y_{n}}\partial f(z)$ satisfying
$\|y_{n}^{*}-x_{0}^{*}\|\le2\varepsilon_{n}/\delta_{n}$; hence
$y_{n}^{*}\to x_{0}^{*}$  in norm and in the weak-star topology.  As in the
proof of \Cref{thm:min-sum}, $y_{n}\to x_{0}$ and $f(y_{n})\to f(x_{0})$.

For each $n$, the membership
\[
y_n^{*}\in s\times w^*-\Limsup_{z\to y_n}\partial f(z)
\]
means that there is a net
\[
(z_{n,i},z_{n,i}^{*})_{i\in I_n}\subset\Graph(\partial f)
\]
such that
\[
z_{n,i}\to y_n,\qquad
z_{n,i}^{*}\to^{\wstar}y_n^{*},\qquad
f(z_{n,i})\to f(y_n).
\]
A standard diagonal-net argument applied to these nets yields a single net
\[
(z_j,z_j^{*})\subset\Graph(\partial f)
\]
such that
\[
z_j\to x_0,\qquad
z_j^{*}\to^{\wstar}x_0^{*},\qquad
f(z_j)\to f(x_0).
\]
For completeness, one may direct the triples consisting of $n$, a strong
neighborhood of $y_n$, a weak-star neighborhood of $y_n^*$, and a positive
function-value tolerance, and then combine this directed set with the filter
$n\to\infty$.
Therefore
\[
x_0^{*}\in s\times w^*-\Limsup_{y\to x_0}\partial f(y).
\]
\end{proof}

If property \textup{[P]} is assumed in the stronger $s\times s^*$-limiting form, then the same proof yields
\begin{equation}
	\eth f(x)\subset s\times s^*-\Limsup_{y\to x}\partial f(y).
\end{equation}
 
For a normed space, write \(\fsub f(x)\) for the Fr\'echet, or regular,
subdifferential:
\[
\fsub f(x):=\left\{x^{*}\in X^{*}\ \middle|\
\liminf_{u\to x,\ u\ne x}
\frac{f(u)-f(x)-\langle u-x,x^{*}\rangle}{\|u-x\|}\ge0\right\},
\]
with the convention \(\fsub f(x)=\emptyset\) when \(f(x)\notin\mathbb R\).  Its
limiting, or Mordukhovich, closure is denoted by
\[
\lsub f(x)\coloneqq s\times w^*-\Limsup_{y\to x}\fsub f(y).
\]

\begin{theorem}[nearby Fr\'echet subgradients]\label{frechet-limiting}
Let \((X,\|\cdot\|)\) be a normed space, let \(f:X\to\Rbar\), and suppose that
\(f(x_0)\in\mathbb R\).
\begin{enumerate}[label=\textup{(\alph*)},leftmargin=2.4em]
\item One has \(\fsub f(x_0)\subset\eth f(x_0)\).
\item If, in addition, \(X\) is Banach, \(f\) is l.s.c. near \(x_0\), and the
Fr\'echet subdifferential is an abstract subdifferential, in the sense of
\textup{[P]}, on a class containing \(f\), then
\[
\eth f(x_0)\subset \lsub f(x_0).
\]
In particular, this applies in the standard Asplund-Banach-space setting for
proper lower semicontinuous functions, where the Fr\'echet subdifferential
satisfies the fuzzy local-minimum principle \textup{[P]} used above; see, for
example, \cite{Mordukhovich2006,Ioffe2017,Kruger2024}.
\end{enumerate}
\end{theorem}

\begin{proof}
For (a), let \(x_0^{*}\in\fsub f(x_0)\).  For each \(n\in\mathbb N\), the
Fr\'echet subgradient inequality with \(\eta=1/n\) gives a radius
\(r_n>0\) such that
\[
f(x_0)+\langle u-x_0,x_0^*\rangle
\le f(u)+\frac1n\|u-x_0\|
\qquad(0<\|u-x_0\|\le r_n).
\]
The same inequality is trivial at \(u=x_0\).  Put
\(\delta_n:=\min\{r_n,1\}\) and
\(\varepsilon_n:=\delta_n/n\).  Then \(\varepsilon_n>0\),
\(\varepsilon_n\to0\), and
\(\varepsilon_n/\delta_n=1/n\to0\).  Moreover, for every
\(u\in\overline B(x_0;\delta_n)\),
\[
f(x_0)+\langle u-x_0,x_0^*\rangle\le f(u)+\varepsilon_n .
\]
By \Cref{equiv-core}, \(x_0^*\in\eth f(x_0)\).

For (b), apply \Cref{thm:min-abstract} with \(\partial=\fsub\).  The Banach
assumption is exactly the completeness hypothesis used in the Ekeland argument
inside \Cref{thm:min-abstract}, and the right-hand side is the definition of
\(\lsub f(x_0)\).  The final Asplund-Banach-space assertion follows from the
standard fuzzy sum rule, or fuzzy minimum principle, for the Fr\'echet
subdifferential on Asplund spaces.
\end{proof}

\begin{theorem}[core subdifferential on Asplund spaces]\label{thm:core-asplund-abstract}
Let \(X\) be an Asplund Banach space, and let \(\mathcal F\) be the class of
proper lower semicontinuous functions \(f:X\to\Rbar\). Then the core
subdifferential \(\eth\) is an abstract subdifferential on \(\mathcal F\) in
the sense of \textup{[P]}; that is, if \(f\in\mathcal F\), \(g:X\to\mathbb R\)
is convex and continuous, and \(x\in X\) is a local minimum point of \(f+g\),
with \(f(x)\in\mathbb R\), then
\[
0\in s\times w^*-\Limsup_{y\to x}\eth f(y)+\fen g(x).
\]
\end{theorem}

\begin{proof}
Since \(X\) is Asplund and \(f\) is proper lower semicontinuous, the
Fr\'echet subdifferential satisfies the fuzzy minimum principle; see, for
instance, Mordukhovich~\cite[Chapter~2]{Mordukhovich2006}.  Hence, because
\(f+g\) has a local minimum at \(x\) and \(g\) is convex continuous,
\[
0\in s\times w^*-\Limsup_{y\to x}\widehat{\partial}f(y)+\fen g(x),
\]
where \(\widehat{\partial}f\) denotes the Fr\'echet subdifferential.  By
\Cref{frechet-limiting}(a), \(\widehat{\partial}f(y)\subset \eth f(y)\) for every
\(y\in X\). Therefore
\[
s\times w^*-\Limsup_{y\to x}\widehat{\partial}f(y)
\subset
s\times w^*-\Limsup_{y\to x}\eth f(y).
\]
Combining the two inclusions gives the asserted \textup{[P]} property for
\(\eth\).
\end{proof}

\begin{remark}[nearby calculus and abstract convexity]\label{rem:nearby-calculus}
Theorem~\ref{frechet-limiting} explains the connection with the nearby-point, or
fuzzy-calculus, philosophy of Fr\'echet subdifferentials: a core slope need not
itself be a Fr\'echet subgradient at \(x\), but in the standard Asplund setting it
is reached as a nearby Fr\'echet subgradient.  This is distinct from the recent
abstract-convexity use of the term abstract subdifferential, where one fixes a
class of abstract affine functions and studies summation, composition, and
abstract monotonicity properties, see for example  \cite{DiazMillanSukhorukovaUgon2025}.  The
present paper instead takes the subdifferential operator as primitive and asks
which slopes are unavoidable under local-minimum and fuzzy nearby-point axioms.
\end{remark}

\eject

\section{Definition and basic properties}

For a metric space $(X,d)$, an extended-real-valued function
$f:X\to\overline{\mathbb{R}}$, a point $x_0\in\operatorname{dom} f$,
and $\varepsilon\ge0$, define the $\varepsilon$-minimality radius of
$f$ at $x_0$ by
\begin{equation}\label{df}
	\delta_f(\varepsilon;x_0)
	:=
	\sup\left\{
	\delta\ge0\;\middle|\;
	f(x_0)\le f(x)+\varepsilon
	\text{ for every }x\in \overline{B}_d(x_0;\delta)
	\right\}\in[0,+\infty].
\end{equation}
For notational convenience, we set, by convention, $	\delta_f(+\infty;x_0)=\infty$. 

\begin{proposition}\label{prop:df}
	Let $(X,d)$ be a metric space, let $f:X\to\Rbar$, and let $x_{0}\in X$ with
	$f(x_{0})\in\mathbb{R}$. Then the following assertions hold.
	\begin{enumerate}[label=\textup{(\roman*)},leftmargin=2.4em]
		
\item Let $\varepsilon\ge0$. If $d(x,x_{0})<\delta_{f}(\varepsilon;x_{0})$, then $f(x_{0})\le f(x)+\varepsilon$. 

Conversely, if $	f(x_{0})\le f(x)+\varepsilon$, for all $x\in X$ with
		$d(x,x_{0})<\delta$, then $\delta\le\delta_{f}(\varepsilon;x_{0})$. 
		
\item The function 		$\varepsilon\mapsto\delta_{f}(\varepsilon;x_{0})$ is nondecreasing and right-continuous on $[0,+\infty)$.
		
\item The function $f$ is l.s.c. at $x_{0}$ if and only if $\delta_{f}(\varepsilon;x_{0})>0$ for every $\varepsilon>0$.
		
\item The point $x_{0}$ is a local minimum point of $f$ if and only if $	\inf_{\varepsilon>0}\delta_{f}(\varepsilon;x_{0})
=\delta_{f}(0;x_{0})>0$.

\item If $f$ is l.s.c. at $x_0$, then, for every
		$\mu\ge0$, $\delta_f$ is upper semicontinuous at $(\mu,x_0)$ along finite-domain points, that is,
		\[
		\limsup_{\varepsilon\to\mu,\,x\to x_0,\, f(x)\in\mathbb R}
		\delta_{f}(\varepsilon;x)
		\le
		\delta_{f}(\mu;x_0).
		\]
	Due to the monotonicity in $\varepsilon$, the case \(x=x_{0}\) is precisely the right-continuity in
	\(\varepsilon\) stated in item~\textup{(ii)}.
		
\item If, for some $C>0$, $g:X\to\mathbb{R}$ satisfies
		\[
		\forall x\in X,\ 
		g(x)\le g(x_{0})+ C d(x,x_{0}),
		\]
	 then,  for every $\varepsilon\ge0$, 

\begin{equation} \label{lls}
	\delta_{f+g}(\varepsilon;x_{0})	\le
\delta_{f}\bigl(\varepsilon+C\delta_{f+g}(\varepsilon;x_{0});x_{0}\bigr).
\end{equation}

	\end{enumerate}
\end{proposition}

\begin{proof}
(i) If $d(x,x_{0})<\delta_{f}(\varepsilon;x_{0})$, choose an admissible
radius $\rho>d(x,x_{0})$ in \eqref{df}; then $x\in \overline{B}_d(x_{0};\rho)$ and
the desired inequality follows.  Conversely, if the inequality holds on the
open ball of radius $\delta$, then every closed ball of radius $\rho<\delta$ is
admissible, so $\rho\le\delta_{f}(\varepsilon;x_{0})$ for all $\rho<\delta$.
Letting $\rho\uparrow\delta$ gives the claim.

(ii) Monotonicity is immediate from the definition.  Put
$$\alpha\coloneqq\lim_{\varepsilon\downarrow\mu}\delta_{f}(\varepsilon;x_{0})=\inf_{\varepsilon>\mu}\delta_{f}(\varepsilon;x_{0}).$$  Since
$\delta_{f}(\mu;x_{0})\le\alpha$, it remains to prove the reverse inequality.
If $r<\alpha$ and $d(x,x_{0})<r$, then by (i)
$f(x_{0})\le f(x)+\varepsilon$ for all $\varepsilon>\mu$.  Letting
$\varepsilon\downarrow\mu$ gives $f(x_{0})\le f(x)+\mu$ on the open ball of
radius $r$, and (i) yields $r\le\delta_{f}(\mu;x_{0})$.  Let
$r\uparrow\alpha$.

(iii) This is exactly the $\varepsilon$--neighborhood characterization of
lower semicontinuity at a finite point, since 
$$
\delta_f(\varepsilon;x_0)>0\ \Leftrightarrow\ \exists r>0,\ 
\forall x\in\overline{B}_d(x_0;r),\ f(x_0)\le f(x)+\varepsilon.
$$

(iv) By (ii),
$\inf_{\varepsilon>0}\delta_{f}(\varepsilon;x_{0})=\delta_{f}(0;x_{0})$.
The condition $\delta_{f}(0;x_{0})>0$ says precisely that the inequality
$f(x_{0})\le f(x)$ holds in a neighborhood of $x_{0}$.

(v) Let
\[
\alpha:=\limsup_{\varepsilon\to\mu,\,x\to x_0,\, f(x)\in\mathbb R}
\delta_{f}(\varepsilon;x).
\]
If $\alpha=0$ then $\alpha\le  \delta_f(\mu;x_0)$.

Suppose now that $\alpha>0$. 
Let $0<r<\alpha$.  Choose $a>0$ with $r+a<\alpha$.  

By the definition of the limsup along finite-domain points, for every $t>0$
and $0<s\le a$ there exist $\varepsilon_{t,s}\ge0$ with
$|\varepsilon_{t,s}-\mu|<t$ and $x_{t,s}$ with $f(x_{t,s})\in\mathbb R$ and
$d(x_{t,s},x_0)<s$ such that
$r+a<\delta_{f}(\varepsilon_{t,s};x_{t,s})$.

Now fix $x\in \overline{B}_d(x_0;r)$. Then  
$$
d(x,x_{t,s})\le d(x,x_0)+d(x_{t,s},x_0)\le r+s\le r+a<\delta_{f}(\varepsilon_{t,s};x_{t,s}).
$$

Hence, by the defining property of $\delta_f$,  $f(x_{t,s})\le f(x)+\varepsilon_{t,s}\le f(x)+\mu+t$. 

Now fix $t>0$ and let $s\downarrow0$. 
Since $x_{t,s}\to x_0$ and $f$ is l.s.c. at $x_0$, 
$$
f(x_0)\le \liminf_{y\to x_0}f(y)\le \liminf_{s\downarrow0}f(x_{t,s})\le f(x)+\mu+t.
$$

Letting $t\downarrow0$, we obtain $f(x_0)\le f(x)+\mu $.

Since  $x\in \overline{B}_d(x_0;r)$ was arbitrary, it follows that 
$$
\forall x\in \overline{B}_d(x_0;r),\ f(x_0)\le f(x)+\mu,
$$
that is, $r$ is admissible in the definition of $\delta_f(\mu;x_0)$. Hence $r\le \delta_f(\mu;x_0)$. 
Since $0<r<\alpha$ was arbitrary, $\alpha\le \delta_f(\mu;x_0)$. 

(vi) Let $\varepsilon\ge0$. If $\delta_{f+g}(\varepsilon;x_{0})\in\{0,+\infty\}$, then the conclusion is immediate. Otherwise, by (i), every $x\in X$ with $d(x,x_{0})<\delta_{f+g}(\varepsilon;x_{0})$ has 
$(f+g)(x_{0})\le(f+g)(x)+\varepsilon$, from which 
\[
f(x_{0})\le f(x)+\varepsilon+g(x)-g(x_{0})
\le f(x)+\varepsilon+C\delta_{f+g}(\varepsilon;x_{0}).
\]
By (i), we find \eqref{lls}.
\end{proof}

\begin{theorem}\label{equiv-core}
Let $(X,s)$ be a metrizable topological vector space and let $d$ be a fixed
compatible metric on $X$.  Let $f:X\to\Rbar$ and $x_{0}\in X$ satisfy
$f(x_{0})\in\mathbb{R}$.  For $x_{0}^{*}\in X^{*}$, the following are
equivalent:
\begin{enumerate}[label=\textup{(\roman*)},leftmargin=2.4em]
\item $x_{0}^{*}\in\eth_d f(x_{0})$.
\item There exist $\varepsilon_{n}>0$ and $\delta_{n}>0$ such that
$\varepsilon_{n}\to0$, $\varepsilon_{n}/\delta_{n}\to0$, and
\[
f(x_{0})+\langle x-x_{0},x_{0}^{*}\rangle\le f(x)+\varepsilon_{n}, 
\]
for every $n$ and every $x\in \overline{B}_d(x_{0};\delta_{n})$.
\item
\[
\limsup_{\varepsilon\downarrow0}
\delta_{f-x_{0}^{*}}(\varepsilon;x_{0})/{\varepsilon}=+
\infty.
\]
\end{enumerate}
\end{theorem}

\begin{proof}
(i) $\Rightarrow$ (ii): choose a sequence $\varepsilon_n\downarrow0$ along which
$\delta_{\varepsilon_n}/\varepsilon_n\to+\infty$ in the definition of the core,
and put $\delta_n:=\delta_{\varepsilon_n}$.

(ii) $\Rightarrow$ (iii) follows from
$\delta_{f-x_0^*}(\varepsilon_n;x_0)\ge\delta_n$.

For (iii) $\Rightarrow$ (i), first note that
$\delta_{f-x_0^*}(\varepsilon;x_0)>0$ for every $\varepsilon>0$: indeed, by
monotonicity in $\varepsilon$, any smaller error at which the radius is positive
also gives a positive radius at $\varepsilon$.  Choose positive numbers
$\delta_\varepsilon<\delta_{f-x_0^*}(\varepsilon;x_0)$ so that
\[
\limsup_{\varepsilon\downarrow0}\delta_\varepsilon/\varepsilon=+\infty .
\]
This is possible because the right-hand side of (iii) is infinite; when
$\delta_{f-x_0^*}(\varepsilon;x_0)=+\infty$ one may choose any finite positive
number below it.  For each $\varepsilon$, the strict inequality
$\delta_\varepsilon<\delta_{f-x_0^*}(\varepsilon;x_0)$ allows us to choose an
admissible radius larger than $\delta_\varepsilon$.  Hence the defining
inequality for $\delta_{f-x_0^*}$ holds on
$\overline B_d(x_0;\delta_\varepsilon)$, that is,
\[
f(x_{0})+\langle x-x_{0},x_{0}^{*}\rangle\le f(x)+\varepsilon
\qquad (x\in\overline B_d(x_0;\delta_\varepsilon)).
\]
This is exactly the definition of $x_0^*\in\eth_d f(x_0)$.
\end{proof}

\begin{definition}[lower scale slope]\label{def:scale-slope}
	Let \((X,d)\) be a metric space and let \(f:X\to\Rbar\) with
	\(f(x)\in\mathbb{R}\). The \emph{lower scale slope} of \(f\) at \(x\) is
	\[
	\mathfrak{s}f(x):=
	\liminf_{r\downarrow0}
	\bigl(f(x)-\inf_{u\in \overline{B}_d(x;r)}f(u)\bigr)^{+}/{r},
	\]
	where \(a^{+}:=\max\{a,0\}\). This quantity records the smallest first-order
	rate, along shrinking metric balls, at which \(f\) can drop below its value at
	\(x\).
\end{definition}

\begin{theorem}[scale-slope characterization]\label{scale-slope}
	Let \((X,s)\) be a metrizable topological vector space and let \(d\) be a fixed
	compatible metric on \(X\). Let \(f:X\to\Rbar\), let \(f(x)\in\mathbb{R}\),
	and let \(x^{*}\in X^{*}\). Put \(F\coloneqq f-x^{*}\). Then
	\begin{equation}\label{ssc}
		x^{*}\in\eth_d f(x)\qquad\Longleftrightarrow\qquad
		\mathfrak{s}F(x)=0.
	\end{equation}
	
	Equivalently, \(x^{*}\in\eth_d f(x)\) if and only if there are arbitrarily small
	radii \(r>0\) on which
	\[
	F(x)-\inf_{u\in \overline{B}_d(x;r)}F(u)=o(r).
	\]
\end{theorem}

\begin{proof}
	Assume first that \(x^{*}\in\eth_d f(x)\). By Theorem~\ref{equiv-core}, there are
	\(\varepsilon_n>0\) and \(\delta_n>0\) such that
	\(\varepsilon_n\to0\), \(\varepsilon_n/\delta_n\to0\), and
	\[
	\forall u\in \overline{B}_d(x;\delta_n),\ 
	F(x)\le F(u)+\varepsilon_n.
	\]
	If a subsequence of \((\delta_n)\) is bounded away from zero, then letting
	\(n\to\infty\) along that subsequence shows that \(x\) is a local minimum point
	of \(F\), because \(\varepsilon_n\to0\). Hence \(\mathfrak{s}F(x)=0\).
	Otherwise, along a subsequence \(\delta_n\downarrow0\),
	\[
	0\le F(x)-\inf_{u\in \overline{B}_d(x;\delta_n)}F(u)\le \varepsilon_n,
	\]
	and division by \(\delta_n\) yields \(\mathfrak{s}F(x)=0\).
	
	Conversely, suppose that \(\mathfrak{s}F(x)=0\). Choose \(r_n\downarrow0\) such
	that
	\[
	\bigl(F(x)-\inf_{u\in \overline{B}_d(x;r_n)}F(u)\bigr)^{+}/r_n\to0.
	\]
	Set
	\[
	\varepsilon_n:=
	\bigl(F(x)-\inf_{u\in \overline{B}_d(x;r_n)}F(u)\bigr)^{+}+r_n/n
	\ge F(x)-\inf_{u\in \overline{B}_d(x;r_n)}F(u).
	\]
	Then \(\varepsilon_n>0\), \(\varepsilon_n/r_n\to0\), and
	\(F(x)\le F(u)+\varepsilon_n\) for every
	\(u\in \overline{B}_d(x;r_n)\). By Theorem~\ref{equiv-core},
	\(x^{*}\in\eth_d f(x)\).
\end{proof}

\begin{corollary}[error-bound reading of the core condition]\label{cor:scale-error}
	Let \((X,s)\) be a metrizable topological vector space and let \(d\) be a fixed
	compatible metric on \(X\). Let \(x^{*}\in\eth_d f(x)\). Then, for every
	\(\eta>0\) and every \(\rho>0\), there is \(r\in(0,\rho)\) such that
	\[
	(f-x^{*})(x)\le \inf_{u\in \overline{B}_d(x;r)}(f-x^{*})(u)+\eta r.
	\]
	Consequently, the affine perturbation \(f-x^{*}\) cannot have a uniform sharp
	local descent gap
	\[
	(f-x^{*})(x)-\inf_{u\in \overline{B}_d(x;r)}(f-x^{*})(u)\ge \alpha r
	\]
	for all sufficiently small \(r>0\) and any \(\alpha>0\).
	
	Thus \(x^{*}\) is a core subgradient precisely when the affine perturbation
	\(f-x^{*}\) has best local descent gap sublinear along arbitrarily small
	scales.
\end{corollary}

\begin{proof}
	Put \(F:=f-x^{*}\). Since \(x^{*}\in\eth_d f(x)\), there exists
	\((\delta_{\varepsilon})_{\varepsilon>0}\subset(0,+\infty)\) such that
$
	\limsup_{\varepsilon\downarrow0}
	\delta_{\varepsilon}/\varepsilon=+\infty
$
	and, for every \(\varepsilon>0\) and every
	\(u\in \overline{B}_d(x;\delta_{\varepsilon})\),
$
	F(x)\le F(u)+\varepsilon
$. 
	
	Let \(\eta>0\) and \(\rho>0\). Choose \(\varepsilon>0\) such that
$\varepsilon<\eta\rho$ and 	$\delta_{\varepsilon}/\varepsilon>\eta^{-1}$.
	Then
	\[
	\varepsilon\eta^{-1}<\min\{\delta_{\varepsilon},\rho\}.
	\]
	Choose
$
	r\in\left(\varepsilon\eta^{-1},\min\{\delta_{\varepsilon},\rho\}\right)$. 
	Then \(r\in(0,\rho)\) and
$
	\overline{B}_d(x;r)\subset \overline{B}_d(x;\delta_{\varepsilon})$. 
	Hence, 	for every \(u\in\overline{B}_d(x;r)\), 
$
	F(x)\le F(u)+\varepsilon 
$. 
Taking the infimum over \(u\in\overline{B}_d(x;r)\), we obtain
	\[
	F(x)\le \inf_{u\in\overline{B}_d(x;r)}F(u)+\varepsilon
	\le \inf_{u\in\overline{B}_d(x;r)}F(u)+\eta r .
	\]
	
	Since \(\eta>0\) and \(\rho>0\) were arbitrary, the quantity
$
	F(x)-\inf_{u\in\overline{B}_d(x;r)}F(u)
$
	is bounded above by \(\eta r\) along arbitrarily small radii \(r>0\). Hence no
	estimate of the form
	\[
	F(x)-\inf_{u\in\overline{B}_d(x;r)}F(u)\ge \alpha r
	\]
	can hold for all sufficiently small \(r>0\) with \(\alpha>0\).
\end{proof}

For the metric \(d\), the \(d\)-core subdifferential can be expressed by using
the \(\varepsilon\)-Fenchel subdifferential, indicator functions of closed
metric balls, and admissible scale families as follows:
\begin{equation}\label{eq:core-vs-fenchel}
	\eth_d f(x)=
	\bigcup_{(\delta_{\varepsilon})\in\mathscr{S}}
	\bigcap_{\varepsilon>0}
	\partial_F^{\varepsilon}\bigl(f+\ind_{\overline{B}_d(x;\delta_{\varepsilon})}\bigr)(x),
\end{equation}
where \(\overline{B}_d(x;\delta_{\varepsilon})\) is the closed \(d\)-ball
centered at \(x\) with radius \(\delta_{\varepsilon}\),
\(\ind_{\overline{B}_d(x;\delta_{\varepsilon})}\) is the indicator function of
this closed ball, and \(\efen h(x)\) denotes the global
\((\varepsilon-)\)Fenchel subdifferential
\[
\efen h(x):=\{x^*\in X^*\mid h(x)+\langle u-x,x^*\rangle\le h(u)(+\varepsilon)
\text{ for all }u\in X\}.
\]
Here \(\mathscr{S}\) is the class
of all positive scale families
\((\delta_{\varepsilon})_{\varepsilon>0}\subset(0,+\infty)\) satisfying
$
\limsup_{\varepsilon\downarrow0}
\delta_{\varepsilon}/{\varepsilon}
=+\infty .
$

\begin{theorem}\label{basic-properties}
	Let \((X,s)\) be a metrizable topological vector space and let \(d\) be a fixed
	compatible metric on \(X\). Let \(f,g:X\to\Rbar\).
	\begin{enumerate}[label=\textup{(\alph*)},leftmargin=2.4em]
		\item For every \(x\in X\), \(\fen f(x)\subset\eth_d f(x)\). 
		
If, in addition, $X$ is normed and \(f\) is convex,
		then \(\eth f(x)=\fen f(x)\) for every \(x\in X\).
		\item If \(x\) is a local minimum point of \(f\), then \(0\in\eth_d f(x)\).
		\item If \(\eth_d f(x)\ne\emptyset\), then \(f\) is l.s.c. at \(x\).
		\item If \(g=f\) near \(x\), then \(\eth_d g(x)=\eth_d f(x)\).
		\item If, in addition, $X$ is normed and \(\eth\) denotes the norm-metric
core, then, for every \(x\in X\), \(\eth f(x)\) is norm closed in \(X^{*}\).
		\item If \(d\) is translation invariant and \(g(x)=f(x+x_{0})\), then \(\eth_d g(0)=\eth_d f(x_{0})\).
		\item For every \(\ell^{*}\in X^{*}\), every \(\lambda>0\), and every
\(\bar x\in X\),
		\[
		\eth_d(\lambda f+\ell^{*})(\bar x)=\lambda\eth_d f(\bar x)+\ell^{*},
		\]
where \(\ell^*\) is identified with the linear function
\(u\mapsto\langle u,\ell^*\rangle\).
	\end{enumerate}
\end{theorem}

\begin{proof}
	(a) If \(x_{0}^{*}\in\fen f(x_{0})\), then
$
	f(x_{0})+\langle x-x_{0},x_{0}^{*}\rangle\le f(x)
$
	for all \(x\in X\); so, for every $\varepsilon>0$ take \(\delta_{\varepsilon}=1\) to get \(x_{0}^{*}\in\eth_d f(x_{0})\).

Assume now that \((X,\|\cdot\|)\) is normed, that \(f\) is convex, and that
\(x_0^*\in\eth f(x_0)\). In particular, \(f(x_0)\in\mathbb R\).
Choose \(\varepsilon_n>0\) and \(\delta_n>0\) as in
Theorem~\ref{equiv-core}, so that \(\varepsilon_n\to0\) and
\(\varepsilon_n/\delta_n\to0\).

Fix \(x\in X\). If \(f(x)=+\infty\), the Fenchel inequality at \(x\) is
trivial. Thus assume that \(f(x)<+\infty\). If \(x=x_0\), there is nothing to
prove. Assume \(x\neq x_0\), and put \(h=x-x_0\). Define
\(t_n:=\min\{\sqrt{\varepsilon_n},\delta_n/\|h\|\}\). Then \(t_n>0\),
\(t_n\to0\), \(x_0+t_nh\in\overline B(x_0;\delta_n)\), and
\(\varepsilon_n/t_n\to0\). Indeed,
\[
\|x_0+t_nh-x_0\|=t_n\|h\|\le \delta_n,
\]
and
\[
\frac{\varepsilon_n}{t_n}
=
\max\left\{\sqrt{\varepsilon_n},\|h\|\frac{\varepsilon_n}{\delta_n}\right\}
\to0.
\]

The core inequality gives
\[
t_n\langle h,x_0^*\rangle
\le
f(x_0+t_nh)-f(x_0)+\varepsilon_n.
\]
Dividing by \(t_n\), we obtain
\[
\langle h,x_0^*\rangle
\le
\frac{f(x_0+t_nh)-f(x_0)}{t_n}
+
\frac{\varepsilon_n}{t_n}.
\]
Letting \(n\to\infty\), and using the existence of the right directional
derivative of the convex function \(f\), we get
\[
\langle x-x_0,x_0^*\rangle
\le
f'_+(x_0;h)\le f(x)-f(x_0).
\]
Thus \(x_0^*\in\fen f(x_0)\).
	
	Parts (b), (d), (f), and (g) follow by rewriting the defining inequality.  For
	(f), translation invariance gives $d(0,u-x_0)=d(x_0,u)$; without this
	hypothesis the assertion can fail.  For (g), divide the error by \(\lambda\)
	in one direction and multiply it by \(\lambda\) in the other.
	
	For (c), if \(x^{*}\in\eth_d f(x)\), then for every \(\varepsilon>0\) there is
	\(\delta_{\varepsilon}>0\) such that, 	for every \(z\in \overline{B}_d(x;\delta_{\varepsilon})\), 
	\[
	f(x)-\varepsilon\le f(z)-\langle z-x,x^{*}\rangle. 
	\]
Since \(x^{*}\) is
	continuous, letting \(z\to x\) gives
	\[
	f(x)-\varepsilon\le \liminf_{z\to x}f(z).
	\]
	Letting \(\varepsilon\downarrow0\), we get that \(f\) is l.s.c. at \(x\).
	
	(e) Assume that \(X\) is normed and that \(d\) is induced by the norm. Fix
	\(x_0\in X\). We prove that \(\eth f(x_0)\) is norm closed in \(X^*\).
	
	Let \((x_n^*)_n\subset \eth f(x_0)\) and suppose that
	\(\|x_n^*-x^*\|\to0\). We show that \(x^*\in\eth f(x_0)\). 
	
	For each
	\(k\in\mathbb N\), choose \(n_k\) such that
	\(\|x_{n_k}^*-x^*\|\le 1/k\). Since \(x_{n_k}^*\in\eth f(x_0)\),
	\Cref{equiv-core} gives admissible pairs \((\tilde\mu_j,\tilde\rho_j)\) with
	\(\tilde\mu_j\to0\) and \(\tilde\mu_j/\tilde\rho_j\to0\).  Choose one pair
	with \(\tilde\mu_j\le 1/k^2\) and
	\(\tilde\mu_j/\tilde\rho_j\le 1/k\), set
	\(\mu_k:=\tilde\mu_j\) and \(\rho_k:=\min\{\tilde\rho_j,1/k\}\).  Then
	\(\rho_k\le 1/k\), \(\mu_k\le1/k^2\), and
	\(\mu_k/\rho_k\le1/k\); moreover, for every
	\(x\in \overline B(x_0;\rho_k)\), 
	$$
	f(x_0)+\langle x-x_0,x_{n_k}^*\rangle\le f(x)+\mu_k.
	$$
	
	Now let \(x\in \overline B(x_0;\rho_k)\). Then \(\|x-x_0\|\le\rho_k\), and
	therefore 
	$$
	\bigl|\langle x-x_0,x^*-x_{n_k}^*\rangle\bigr|
	\le \|x-x_0\|\,\|x^*-x_{n_k}^*\|\le \rho_k/k.
	$$
	
	 Hence
	$$
	f(x_0)+\langle x-x_0,x^*\rangle
	=f(x_0)+\langle x-x_0,x_{n_k}^*\rangle+\langle x-x_0,x^*-x_{n_k}^*\rangle
	\le f(x)+\mu_k+\rho_k/k.
	$$
	
	Set \(\beta_k:=\mu_k+\rho_k/k\). Then \(\beta_k\to0\) and
	\(\beta_k/\rho_k=\mu_k/\rho_k+1/k\le2/k\to0\). Thus, for every
	\(x\in \overline B(x_0;\rho_k)\), 
	$$
		f(x_0)+\langle x-x_0,x^*\rangle\le f(x)+\beta_k, 
	$$
with \(\beta_k\to0\) and
	\(\beta_k/\rho_k\to0\). By Theorem~\ref{equiv-core},
	\(x^*\in\eth f(x_0)\). Hence \(\eth f(x_0)\) is norm closed in \(X^*\).
\end{proof}

\begin{remark}\label{rem:indicator-example}
The core subdifferential can differ from the Fenchel subdifferential.  Let
$(X,\|\cdot\|)$ be a normed space and
$$S:=\{x\in X\mid\|x\|=1\}\cup\{0\}.$$  At $0$,
$\fen\ind_{S}(0)=\{0\}$, since for every $u\in X$ with $\|u\|=1$, $u,-u\in S$, hence any $x^*\in \fen\ind_S(0)$ satisfies $\langle u,x_0^*\rangle\le0$ and $\langle -u,x_0^*\rangle\le0$, so $x^*=0$.
On the other hand, $S\cap B(0;1/2)=\{0\}$, so every $x_{0}^{*}\in X^{*}$ satisfies the core
inequality on $\overline B(0;1/2)$ for every $\varepsilon>0$. Thus $\eth\ind_{S}(0)=X^{*}$. 
\end{remark}

\begin{remark}\label{minus-norm}
For every normed space,
\begin{equation}\label{eq:minus-norm-empty}
\eth(-\|\cdot\|)(0)=\emptyset.
\end{equation}
Indeed, if $x_{0}^{*}\in\eth(-\|\cdot\|)(0)$, then for every $\varepsilon>0$ there is $\delta_\varepsilon>0$ such that, for $\|x\|\le\delta_\varepsilon$, $\langle x,x_0^*\rangle\le -\|x\|+\varepsilon$. For all $\|v\|=1$, take $x=\delta_\varepsilon v$ to get 
\[
1+\langle v,x_{0}^{*}\rangle\le\varepsilon/{\delta_{\varepsilon}}.
\]
Taking the supremum over $\|v\|=1$ and letting $\varepsilon\downarrow0$ along
a subnet with $\varepsilon/\delta_{\varepsilon}\to0$ gives
$1+\|x_{0}^{*}\|\le0$, a contradiction. Therefore $\eth(-\|\cdot\|)(0)=\emptyset$.

Consequently, the exact sum rule
$\eth(f+g)(x)\subset\eth f(x)+\eth g(x)$ cannot hold for arbitrary Lipschitz $f$ and
continuous convex $g$. Indeed, take $f=-\|\cdot\|$, $g=\|\cdot\|$, and
$x=0$. Then $f+g=0$, $\eth (f+g)(0)={0}$; but $\eth f(0)=\emptyset$, so $\eth f(0)+\eth g(0)=\emptyset$. 
\end{remark}

\begin{proposition}\label{GaFr}
Let $(X,\|\cdot\|)$ be a normed space, let $x_{0}\in X$, and let
$g:X\to\Rbar$.
\begin{enumerate}[label=\textup{(\alph*)},leftmargin=2.4em]
\item If $x_{0}$ is a local maximum point of $g$, then
$\eth g(x_{0})\subset\{0\}$.
\item If $g$ is G\^ateaux differentiable at $x_{0}$, then
$\eth g(x_{0})\subset\{\nabla g(x_{0})\}$.
\item If $g$ is Fr\'echet differentiable at $x_{0}$, then
$\eth g(x_{0})=\{\nabla g(x_{0})\}$.
\item If $g$ is proper convex and $\eth(-g)(x_{0})\ne\emptyset$, then $g$ is
Fr\'echet differentiable at $x_{0}$ and
\[
\eth g(x_{0})=-\eth(-g)(x_{0})=\{\nabla g(x_{0})\}.
\]
\end{enumerate}
\end{proposition}

\begin{proof}
(a) If $x_{0}^{*}\in\eth g(x_{0})$, choose
$\varepsilon_{n},\delta_{n}$ as in Theorem~\ref{equiv-core}.  
Let $\delta_{0}>0$ be such that $x_0$ is a maximum point of $g$ on $\overline{B}(x_0;\delta_0)$. Let $t_n\coloneqq\min\{\sqrt{\varepsilon_n},\delta_n\}$. Then $t_n\to0$ and $\varepsilon_n/t_n=\max\{\sqrt{\varepsilon_n},\varepsilon_n/\delta_n\}\to0$. For $n$ large, $t_n\le\delta_0$. Also $t_n\le\delta_n$.  The defining inequality used for $x=x_0\pm t_n v\in\overline{B}(x_0;\delta_n)$, where $\|v\|=1$, and the local maximality of $x_0$ provide 
$\pm t_{n}\langle v,x_{0}^{*}\rangle\le g(x_0\pm t_n v)-g(x_0)+
\varepsilon_{n}\le \varepsilon_{n}$.  Hence
$|\langle v,x_{0}^{*}\rangle|\le\varepsilon_{n}/t_{n}$ for every $\|v\|=1$.
Taking the supremum over $\|v\|=1$ gives
$\|x_{0}^{*}\|\le\varepsilon_{n}/t_{n}\to0$. Therefore $x_{0}^{*}=0$.

(b)  Assume that $g$ is G\^ateaux differentiable at $x_{0}$ and let $x_0^*\in\eth g(x_0)$. Choose
$\varepsilon_{n},\delta_{n}$ as in Theorem~\ref{equiv-core}. 

If $\inf_{n}\delta_{n}>0$, then $x_{0}$ is a local minimum point of
$g-x_{0}^{*}$, and the G\^ateaux derivative of this function at $x_{0}$ is zero.  Hence $x_0^*=\nabla g(x_{0})$.

Otherwise pass to a subsequence with $\delta_{n}\to0$.  Fix
$v\in X$ and put $a_{v}:=\max\{1,\|v\|\}$ and
$t_{n}:=\delta_{n}/a_{v}$.  Then $\|t_{n}v\|\le\delta_{n}$, $t_{n}\to0$,
and the core inequality gives
\[
\langle v,x_{0}^{*}\rangle
\le \frac{g(x_{0}+t_{n}v)-g(x_{0})}{t_{n}}
+\frac{\varepsilon_{n}}{t_{n}}.
\]
Since $\varepsilon_{n}/t_{n}=a_{v}\varepsilon_{n}/\delta_{n}\to0$, and the
same argument applied to $-v$ gives the reverse inequality, the G\^ateaux
derivative satisfies $\langle v,x_{0}^{*}\rangle=\langle v,\nabla g(x_{0})\rangle$ for every $v\in X$.  Hence
$x_{0}^{*}=\nabla g(x_{0})$.

(c) Assume that  $g$ is Fr\'echet differentiable at $x_{0}$. 
Inclusion ``$\subset$'' follows from (b).  

For the reverse inclusion,
let $a:=\nabla g(x_{0})$.  For each $k$ choose $0<\rho_{k}\le1$ such that
\[
|g(x)-g(x_{0})-\langle x-x_{0},a\rangle|
\le k^{-2}\|x-x_{0}\|\qquad(\|x-x_{0}\|\le\rho_{k}).
\]
Put $\varepsilon_{k}:=k^{-2}\rho_{k}$.  Then
$\varepsilon_{k}\to0$ and $\varepsilon_{k}/\rho_{k}=k^{-2}\to0$, which proves
$a\in\eth g(x_{0})$.

(d) Let $x_{0}^{*}\in\eth(-g)(x_{0})$ and put
$h(u):=g(u)+\langle u,x_{0}^{*}\rangle$.  Then $h$ is proper convex and, for
sequences $\varepsilon_{n},\delta_{n}$ with
$\varepsilon_{n}/\delta_{n}\to0$,
\[
h(x)\le h(x_{0})+\varepsilon_{n}\qquad(\|x-x_{0}\|\le\delta_{n}).
\]

Now fix \(u\in X\) with \(0<\|u\|\le \delta_n\). Put
\(v:=u/\|u\|\) and \(t:=\|u\|/\delta_n\). Then \(\|v\|=1\),
\(0<t\le1\), and \(x_0+u=(1-t)x_0+t(x_0+\delta_n v)\). By convexity of
\(h\),
\[
h(x_0+u)
\le
(1-t)h(x_0)+t h(x_0+\delta_n v).
\]
Since \(x_0+\delta_n v\in \overline B(x_0;\delta_n)\), the upper estimate
gives
$
h(x_0+\delta_n v)\le h(x_0)+\varepsilon_n
$. 
Hence
\[
h(x_0+u)-h(x_0)
\le
t\varepsilon_n
=
\frac{\varepsilon_n}{\delta_n}\|u\|.
\]
For the lower bound, use
\(x_0=\frac12(x_0+u)+\frac12(x_0-u)\). By convexity,
$
h(x_0)
\le
\frac12 h(x_0+u)+\frac12 h(x_0-u)
$. 
Therefore
\[
h(x_0+u)-h(x_0)
\ge
-\bigl(h(x_0-u)-h(x_0)\bigr).
\]
Applying the previous upper bound to \(-u\), we get
\[
h(x_0-u)-h(x_0)
\le
\frac{\varepsilon_n}{\delta_n}\|u\|.
\]
Thus
\[
h(x_0+u)-h(x_0)
\ge
-\frac{\varepsilon_n}{\delta_n}\|u\|.
\]
Combining the two estimates,
\[
|h(x_0+u)-h(x_0)|
\le
\frac{\varepsilon_n}{\delta_n}\|u\|
\qquad
(0<\|u\|\le\delta_n).
\]
Since \(\varepsilon_n/\delta_n\to0\), it follows that \(h\) is
Fr\'echet differentiable at \(x_0\) with derivative \(0\). Consequently,
\(g\) is Fr\'echet differentiable at \(x_0\), with
\(\nabla g(x_0)=-x_0^*\).

  Part (c) gives
$\eth g(x_{0})=\{\nabla g(x_{0})\}$ and
$\eth(-g)(x_{0})=\{-\nabla g(x_{0})\}$.
\end{proof}

\begin{lemma}\label{lem:local-min}
Let $(X,\|\cdot\|)$ be a normed space, let $f:X\to\Rbar$, and let
$g:X\to\mathbb{R}\cup\{+\infty\}$.  If $x_{0}$ is a local minimum point of
$f+g$ and $f(x_{0})\in\mathbb{R}$, then
\[
\eth(-g)(x_{0})\subset\eth f(x_{0}).
\]
\end{lemma}

\begin{proof}
Let $x_{0}^{*}\in\eth(-g)(x_{0})$.  Then $g(x_{0})\in\mathbb{R}$ and there
are $\varepsilon_{n}>0$, $\delta_{n}>0$ such that
$\varepsilon_{n}\to0$, $\varepsilon_{n}/\delta_{n}\to0$, and
\[
-g(x_{0})\le -g(x)-\langle x-x_{0},x_{0}^{*}\rangle+\varepsilon_{n}
\]
for $\|x-x_{0}\|\le\delta_{n}$.  
In particular, $g(x)<+\infty$ for such $x$, since 
otherwise the right-hand side would be $-\infty$. 

If $x_{0}$ is a minimum point of $f+g$ on
$B(x_{0};\delta_{0})$, where $\delta_0>0$, then for
$r_{n}:=\min\{\delta_{0},\delta_{n}\}>0$ and $\|x-x_{0}\|\le r_{n}$, we have $f(x_0)+g(x_0)\le f(x)+g(x)$ and also 
\[
f(x_{0})+\langle x-x_{0},x_{0}^{*}\rangle\le f(x)+\varepsilon_{n}.
\]
Since $r_{n}/\varepsilon_{n}\to+
\infty$, this proves $x_{0}^{*}\in\eth f(x_{0})$.
\end{proof}

\begin{theorem}\label{thm:sum-rule}
Let $(X,\|\cdot\|)$ be a normed space and let $f,g:X\to\Rbar$.  If $x_{0}$ is
a local minimum point of $f+g$, $f(x_{0})\in\mathbb{R}$, and $g$ is
Fr\'echet differentiable at $x_{0}$, then
\[
-\nabla g(x_{0})\in\eth f(x_{0}).
\]
In particular, if $g$ is proper convex and $\eth(-g)(x_{0})\ne\emptyset$, then
\[
0\in\eth f(x_{0})+\eth g(x_{0}).
\]
\end{theorem}

\begin{proof}
Since $g$ is Fr\'echet differentiable at $x_{0}$,
$-\nabla g(x_{0})\in\eth(-g)(x_{0})$ by Theorem~\ref{GaFr}(c) applied to $-g$.
The first assertion follows from Lemma~\ref{lem:local-min}.  The second assertion
follows from Theorem~\ref{GaFr}(d), which gives
$\eth g(x_{0})=\{\nabla g(x_{0})\}$.
\end{proof}

The same abstract-minimality principle gives the comparison with the
Clarke--Rockafellar subdifferential.

\begin{theorem}\label{thm:Clarke}
Let $(X,\|\cdot\|)$ be a Banach space and let $f:X\to\Rbar$ be l.s.c. near
$x\in X$.  Then
\[
\eth f(x)\subset\partial_{C}f(x),
\]
where $\partial_{C}f$ denotes the Clarke--Rockafellar subdifferential for
extended-real lower semicontinuous functions.
\end{theorem}
\begin{proof}
If $f(x)\notin\mathbb{R}$, then $\eth f(x)=\emptyset$.  Assume that
$f(x)\in\mathbb R$.  Two standard facts about the Clarke--Rockafellar
subdifferential are being used here.  First, it satisfies the abstract fuzzy
minimum principle \textup{[P]} on Banach spaces for lower semicontinuous
functions; this follows from the Clarke--Rockafellar calculus for sums with
continuous convex functions.  Second, its graph is closed for
strong $\times$ weak-star convergence when the base points converge strongly and
the function values converge.  See, for example,
\cite{MR709590,MR1921556,RockafellarWets1998}.

By the first fact and \Cref{thm:min-abstract}, applied with
$\partial=\partial_C$,
\[
\eth f(x)\subset s\times w^*-\Limsup_{y\to x}\partial_C f(y).
\]
Thus any $x^*$ in the right-hand side is represented by a net
$(y_i,y_i^*)\subset\Graph(\partial_C f)$ with
$y_i\to x$, $f(y_i)\to f(x)$, and $y_i^*\to^{w^*}x^*$.  The closed-graph
property gives $x^*\in\partial_C f(x)$, and the inclusion follows.
\end{proof}
For instance,
\[
\eth(-\|\cdot\|)(0)=\emptyset\subsetneqq\partial_{C}(-\|\cdot\|)(0)=B_{X^{*}},
\]
where $B_{X^{*}}$ is the closed dual unit ball.  Thus, even for locally
Lipschitz functions, the core subdifferential can be strictly smaller than the
Clarke--Rockafellar subdifferential. 

\section{Comparisons, examples, and applications}\label{sec:comparisons-applications}

The preceding results show that \(\eth\) is a small but robust object among
common subdifferentials.  This section records concrete consequences that help
locate it in the nonsmooth-analysis landscape.

\subsection{Standard generalized gradients}

The preceding results place \(\eth\) between exact regular support and the
standard limiting or convexified generalized gradients.  At the regular level,
\Cref{frechet-limiting} gives
\[
\fsub f(x)\subset\eth f(x)
\]
whenever \(f(x)\) is finite.  Thus every Fr\'echet, or regular, subgradient is a
core subgradient.  This inclusion does not require completeness or any fuzzy
calculus: a Fr\'echet support is already an affine lower support up to an
\(o(\|u-x\|)\) error, and therefore up to an \(o(\delta)\) error on balls of
radius \(\delta\).

The converse comparison is not pointwise with \(\fsub\), but nearby.  By
\Cref{thm:min-abstract}, if \(\partial\) is an abstract subdifferential
satisfying the local-minimum principle \textup{[P]}, then, for every \(f\) in
its class which is lower semicontinuous near \(x\),
\[
\eth f(x)\subset s\times w^*-\Limsup_{y\to x}\partial f(y).
\]
In this sense \(\eth\) is the minimal local-support object forced by all
abstract subdifferentials satisfying the standard fuzzy nearby-point minimum
principle.

Applying this abstract-minimality result to the Fr\'echet subdifferential gives,
whenever the Fr\'echet fuzzy calculus applies,
\[
\eth f(x)\subset\lsub f(x),
\]
where
\[
\lsub f(x):=s\times w^*-\Limsup_{y\to x}\fsub f(y)
\]
is the limiting, or Mordukhovich, subdifferential.  In particular, this inclusion
holds for proper lower semicontinuous functions on Asplund spaces.  Moreover,
\Cref{thm:core-asplund-abstract} shows that, on Asplund Banach spaces, the core
subdifferential itself satisfies the same abstract nearby-point principle on the
class of proper lower semicontinuous functions.

The same principle also yields the Clarke--Rockafellar comparison.  If \(X\) is
a Banach space and \(f\) is lower semicontinuous near \(x\), then
\Cref{thm:Clarke} gives
\[
\eth f(x)\subset\csub f(x),
\]
where \(\csub\) denotes the Clarke--Rockafellar subdifferential.  Thus core
subgradients are also Clarke--Rockafellar subgradients, although the converse
may fail.

Finally, in finite-dimensional locally Lipschitz problems, the comparison with
Goldstein subdifferentials follows from the Clarke inclusion.  Namely,
\Cref{prop:goldstein} below gives, for every \(r>0\),
\[
\eth f(x)\subset\csub f(x)\subset\gsub^{r}f(x).
\]
Hence core stationarity implies Goldstein stationarity at every radius, but not
conversely.

Consequently, in the standard Banach settings where the fuzzy calculus is
available, the core subdifferential sits between the exact Fr\'echet
subdifferential and the usual limiting or convexified generalized gradients.  It
is always contained in the Clarke--Rockafellar subdifferential for lower
semicontinuous Banach-space functions, and in finite-dimensional locally
Lipschitz problems it is contained in every Goldstein subdifferential.  Its
robustness is therefore not a convexification or relaxation property.  Rather,
it is a nearby-point robustness: core slopes are affine lower supports stable on
scales whose supporting radius dominates the error, and every abstract
subdifferential satisfying the standard fuzzy minimum principle must recover
them as nearby limiting subgradients.  The examples below show that the
inclusions can be strict; in particular, \(\eth\) may exclude generalized slopes
arising from sharp local maxima or from relaxation alone.

\begin{example}[finite-dimensional test cases]\label{ex:finite-dimensional-tests}
	Let \(X=\mathbb R^{m}\) with the Euclidean norm.
	\begin{enumerate}[label=\textup{(\alph*)},leftmargin=2.4em]
		
		\item If
		\[
		f(x)=\max_{1\le i\le N}\{a_i+\langle u_i,x\rangle\},
		\]
		then \(f\) is convex and
		\[
		\eth f(x)=\fen f(x)
		=\co\{u_i\mid a_i+\langle u_i,x\rangle=f(x)\}.
		\]
		Thus the core gives the usual active-slope formula for convex piecewise-linear
		functions.
		
		\item For \(f(t)=-|t|\) on \(\mathbb R\),
		\[
		\eth f(0)=\emptyset,
		\qquad
		\csub f(0)=[-1,1].
		\]
		Hence Clarke and Goldstein stationarity may hold at a sharp local maximum where
		the core detects no lower support.
		
		\item More generally, for \(f_{\alpha}(t)=-|t|^{\alpha}\) at \(0\), with
		\(\alpha>0\),
		\[
		\eth f_{\alpha}(0)=
		\begin{cases}
			\{0\},&\alpha>1,\\
			\emptyset,&0<\alpha\le1.
		\end{cases}
		\]
		The threshold \(\alpha=1\) reflects exactly whether the negative cusp can be
		controlled by an error \(o(\delta)\) on a radius-\(\delta\) ball.
		
		\item If \(C\subset X\) is closed and convex, then the distance function
		\[
		d_C(x):=\dist(x,C)
		\]
		is convex, and hence
		\[
		\eth d_C(x)=\fen d_C(x).
		\]
		In particular, if \(x\in C\), then
		\[
		\eth d_C(x)=N_C(x)\cap B_{X^*},
		\]
		where \(N_C\) is the usual convex normal cone.
		
		\item If
		\[
		S=\{0\}\cup\{x\in X\mid \|x\|=1\},
		\]
		then
		\[
		\eth\ind_S(0)=X^*.
		\]
		This illustrates that isolated feasible points have maximal core normal cone,
		just as they have no nearby feasible first-order geometry to restrict affine
		supports.
		
	\end{enumerate}
\end{example}

\begin{proof}[Verification of Example~\ref{ex:finite-dimensional-tests}(c)]
Let $a\in\mathbb R$.  If $a\in\eth f_\alpha(0)$, then there are
$\varepsilon_n>0$ and $\delta_n>0$ with $\varepsilon_n\to0$,
$\varepsilon_n/\delta_n\to0$, and
\[
a t\le -|t|^\alpha+\varepsilon_n
\qquad (|t|\le\delta_n).
\]
If $a>0$, taking $t=\delta_n$ gives
\[
a\delta_n+\delta_n^\alpha\le\varepsilon_n,
\]
and if $a<0$, taking $t=-\delta_n$ gives the same estimate with $|a|$ in
place of $a$.  After division by $\delta_n$ this contradicts
$\varepsilon_n/\delta_n\to0$.  Hence every core subgradient must
satisfy $a=0$.

For $a=0$, the core inequality is
\[
|t|^\alpha\le\varepsilon_n
\qquad (|t|\le\delta_n).
\]
In particular, $\delta_n^\alpha\le\varepsilon_n$.  If
$0<\alpha\le1$, then $\delta_n\to0$ and
\[
\delta_n^{\alpha-1}\le\varepsilon_n/\delta_n\to0,
\]
which is impossible: the left side is identically $1$ when $\alpha=1$ and
diverges when $0<\alpha<1$.  Thus $0\notin\eth f_\alpha(0)$ for
$0<\alpha\le1$.

If $\alpha>1$, choose any $\delta_n\downarrow0$ and set
$\varepsilon_n:=\delta_n^\alpha$.  Then
$\varepsilon_n/\delta_n=\delta_n^{\alpha-1}\to0$ and
$|t|^\alpha\le\varepsilon_n$ whenever $|t|\le\delta_n$.
Hence $0\in\eth f_\alpha(0)$ for $\alpha>1$.  This proves the
claimed formula.
\end{proof}

\subsection{Constrained problems and variational inequalities}

For a set \(C\subset X\), define the \emph{core normal cone} by
\[
N_C^{\eth}(x):=\eth\ind_C(x).
\]
This notation is useful because the differentiable core sum rule immediately
produces first-order necessary conditions for constrained problems.

\begin{proposition}[core normal optimality conditions]\label{prop:core-normal-apps}
	Let \((X,\|\cdot\|)\) be a normed space and let \(C\subset X\).
	\begin{enumerate}[label=\textup{(\alph*)},leftmargin=2.4em]
		
		\item If \(\varphi:X\to\mathbb R\) is Fr\'echet differentiable at
		\(x_0\in C\) and \(x_0\) is a local minimizer of \(\varphi\) over \(C\), then
		\[
		-\nabla\varphi(x_0)\in N_C^{\eth}(x_0).
		\]
		
		\item If \(F:C\to X^*\) and \(x_0\in C\) satisfies the local variational
		inequality
		\[
		\langle y-x_0,F(x_0)\rangle\ge0
		\]
		for all \(y\in C\) near \(x_0\), then
		\[
		-F(x_0)\in N_C^{\eth}(x_0).
		\]
		
		\item If \(C\) is convex, then
		\[
		N_C^{\eth}(x)=N_C(x)
		\]
		for every \(x\in C\), where \(N_C(x)\) is the usual convex normal cone.
		If \(x\notin C\), then \(N_C^{\eth}(x)=\emptyset\).
		
	\end{enumerate}
\end{proposition}

\begin{proof}
	For (a), the assumption says that \(x_0\) is a local minimum point of
	\(\ind_C+\varphi\).  Since \(\ind_C(x_0)=0\) and \(\varphi\) is Fr\'echet
	differentiable at \(x_0\), \Cref{thm:sum-rule}, applied with
	\(f=\ind_C\) and \(g=\varphi\), gives
	\[
	-\nabla\varphi(x_0)\in\eth\ind_C(x_0)=N_C^{\eth}(x_0).
	\]
	
	For (b), define the continuous linear function
	\[
	\ell(y):=\langle y,F(x_0)\rangle .
	\]
	The displayed variational inequality says exactly that \(x_0\) is a local
	minimizer of \(\ell\) over \(C\).  Applying part (a) to \(\varphi=\ell\), and
	using \(\nabla\ell=F(x_0)\), gives
	\[
	-F(x_0)\in N_C^{\eth}(x_0).
	\]
	
	For (c), if \(x\in C\), then \(\ind_C(x)=0\) and, since \(C\) is convex,
	\(\ind_C\) is convex.  Hence \Cref{basic-properties} gives
	\[
	\eth\ind_C(x)=\fen\ind_C(x).
	\]
	But the Fenchel subdifferential of \(\ind_C\) is precisely the convex normal
	cone:
	\[
	\fen\ind_C(x)=N_C(x).
	\]
	Therefore \(N_C^{\eth}(x)=N_C(x)\) for \(x\in C\).  If \(x\notin C\), then
	\(\ind_C(x)=+\infty\), and by definition \(\eth\ind_C(x)=\emptyset\).
\end{proof}

\subsection{A scale-wise level-set error estimate}

The definition of \(\eth\) is pointwise in the slope.  For error-bound purposes,
a uniform version is more useful.

\begin{definition}[uniform core ball]\label{def:uniform-core-ball}
	Let \((X,\|\cdot\|)\) be a normed space, let \(f:X\to\Rbar\), let
	\(\bar x\in X\) with \(f(\bar x)\in\mathbb R\), and let \(\alpha>0\).  We say
	that \(f\) has a \emph{uniform \(\alpha\)-core ball} at \(\bar x\) if there are
	sequences \(\varepsilon_n>0\) and \(\delta_n>0\) such that
	\[
	\varepsilon_n\to0,
	\qquad
	\frac{\varepsilon_n}{\delta_n}\to0,
	\]
	and
	\[
	f(\bar x)+\langle x-\bar x,u^*\rangle\le f(x)+\varepsilon_n
	\]
	for every \(u^*\in\alpha B_{X^*}\), every \(n\), and every
	\(x\in B(\bar x;\delta_n)\).
\end{definition}

\begin{proposition}[level-set localization]\label{level-error}
	If \(f\) has a uniform \(\alpha\)-core ball at \(\bar x\), then, for every
	\(n\) and every \(t\ge0\),
	\[
	\{x\in B(\bar x;\delta_n)\mid f(x)\le f(\bar x)+t\}
	\subset
	B\left(\bar x;\frac{t+\varepsilon_n}{\alpha}\right).
	\]
	Thus lower level sets are localized around \(\bar x\), with an error term
	\(\varepsilon_n=o(\delta_n)\) on the scale \(\delta_n\).
\end{proposition}

\begin{proof}
	Let \(x\in B(\bar x;\delta_n)\) and set \(h:=x-\bar x\).  If \(h=0\), there is
	nothing to prove.  If \(h\ne0\), the Hahn--Banach theorem gives
	\(\xi^*\in B_{X^*}\) such that
	\[
	\langle h,\xi^*\rangle=\|h\|.
	\]
	Applying the uniform core inequality with \(u^*=\alpha\xi^*\) gives
	\[
	f(\bar x)+\alpha\|h\|\le f(x)+\varepsilon_n.
	\]
	If also \(f(x)\le f(\bar x)+t\), then
	\[
	\alpha\|h\|\le t+\varepsilon_n.
	\]
	Hence
	\[
	\|x-\bar x\|\le \frac{t+\varepsilon_n}{\alpha},
	\]
	which proves the asserted inclusion.
\end{proof}

\subsection{Goldstein stationarity}

In finite-dimensional locally Lipschitz optimization, Goldstein's
\(r\)-subdifferential is commonly written as
\[
\gsub^r f(x):=
\cl\,\co\bigcup_{\|y-x\|\le r}\csub f(y),
\qquad r>0.
\]
It is a relaxed Clarke-type object designed for algorithmic stationarity
certificates
\cite{Goldstein1977,BurkeLewisOverton2005,ZhangLinJegelkaSraJadbabaie2020,LinZhengJordan2022,Gebken2025}.

\begin{proposition}[core stationarity is stronger than Goldstein stationarity]\label{prop:goldstein}
	Let \(X=\mathbb R^m\), let \(f:X\to\mathbb R\) be locally Lipschitz, and let
	\(r>0\).  Then
	\[
	\eth f(x)\subset\csub f(x)\subset\gsub^r f(x).
	\]
	Consequently, \(0\in\eth f(x)\) implies \(0\in\gsub^r f(x)\) for every
	\(r>0\).  The converse fails even in one dimension.
\end{proposition}

\begin{proof}
	The first inclusion is \Cref{thm:Clarke}.  The second follows directly from the
	definition of \(\gsub^r f(x)\), since the union defining \(\gsub^r f(x)\)
	contains \(\csub f(x)\), corresponding to the point \(y=x\).
	
	For failure of the converse, take \(f(t)=-|t|\) at \(0\).  Then
	\[
	0\in\csub f(0)\subset\gsub^r f(0)
	\]
	for every \(r>0\), while
	\[
	\eth f(0)=\emptyset
	\]
	by \Cref{minus-norm}.  Thus Goldstein stationarity does not imply core
	stationarity.
\end{proof}

\section{Conclusion and further directions}\label{sec:conclusion}

The core subdifferential is best viewed as a minimal robust-support object.  It
is not intended to replace Fr\'echet, limiting, Clarke--Rockafellar, or Goldstein
subdifferentials.  Rather, it identifies the slopes that every subdifferential
with a reasonable local-minimum or nearby-point principle must recover.  This is
why the main results are inclusions into graph closures of abstract
subdifferentials, while the examples show that many Clarke or Goldstein slopes
are deliberately excluded.

The most promising further developments are the following.  First, one can ask
for sharper sufficient conditions under which the pointwise core contains a
uniform core ball, because \Cref{scale-slope}, \Cref{level-error} then connect
sublinear scale descent with level-set error bounds.  Second, the core normal
cone \(N_C^{\eth}\) can be used to formulate local variational inequalities and
constrained optimality conditions before passing to larger limiting or Clarke
normal cones.  Third, the comparison with Goldstein stationarity suggests an
algorithmic diagnostic: a point may be Goldstein-stationary because of
convexification over a neighborhood, while an empty or small core
subdifferential records the absence of a robust lower affine support at the base
point.  This distinction may be useful when separating minimizing stationarity
from maximizer-type nonsmooth criticality.

\end{document}